\documentclass[12pt, reqno, a4paper, final]{amsart}

\usepackage{amsmath}                     
\usepackage{a4}
\usepackage{amssymb}                     
\usepackage{amsthm}
\usepackage{hyperref}
\usepackage{marvosym}                   
\usepackage{mathrsfs}									
\usepackage[all]{xy}										
\usepackage{url}												
\usepackage{verbatim,epsfig}

\theoremstyle{plain}

\newtheorem{thm}{Theorem}[section] 
\newtheorem{cor}[thm]{Corollary}
\newtheorem{lem}[thm]{Lemma}
\newtheorem{defi}[thm]{Definition}
\newtheorem{prop}[thm]{Proposition}
\newtheorem*{thm*}{Theorem}
\newtheorem*{prop*}{Proposition}
\newtheorem*{cor*}{Corollary}

\theoremstyle{definition}

\newtheorem{ex}[thm]{Example}
\newtheorem{notation}[thm]{Notation}

\newtheorem{rem}[thm]{Remark}

\newcommand{\NN}{{\mathbb N}}
\newcommand{\ZZ}{{\mathbb Z}}

\newcommand{\RR}{{\mathbb R}}
\newcommand{\CC}{{\mathbb C}}

\newcommand{\GG}{{\mathbb G}}
\newcommand{\HH}{{\mathbb H}}

\renewcommand{\L}{{\mathscr L}}

\newcommand{\K}{{\mathcal K}}

\newcommand{\U}{{\mathcal U}}

\renewcommand{\max}{{\operatorname{max}}}

\newcommand{\ip}[2]{\left\langle {#1} \mid {#2}\right \rangle}

\newcommand{\ipp}{\ip{\cdot}{\cdot}}

\newcommand{\varps}{{\varepsilon}}

\newcommand{\htens}{\bar{\otimes}}

\newcommand{\tens}{\otimes}

\newcommand{\To}{\longrightarrow}

\newcommand{\red}{{\operatorname{r}}}
\newcommand{\Tor}{\operatorname{Tor}}
\newcommand{\op}{{\operatorname{{op}}}}

\newcommand{\Hom}{\operatorname{Hom}}

\newcommand{\id}{\operatorname{id}}

\newcommand{\bet}{\beta^{(2)}}

\newcommand{\Irred}{\operatorname{Irred}}

\renewcommand{\int}{{\operatorname{int}}}

\newcommand{\Pol}{{\operatorname{Pol}}}

\newcommand{\I}{{\operatorname{I}}}
\newcommand{\Ext}{{\operatorname{Ext}}}
\renewcommand{\max}{{\operatorname{u}}}

\newcommand{\maxx}{{\operatorname{max}}}

\newcommand{\one}{{(1)}}
\newcommand{\two}{{(2)}}
\newcommand{\Corep}{{\operatorname{Corep}}}
\newcommand{\T}{{\operatorname{(T)}}}

\renewcommand{\leq}{\leqslant}
\renewcommand{\geq}{\geqslant}

\title{A cohomological description of property (T) for quantum groups}

\author{David Kyed} 
\address{David Kyed,
Mathematisches Institut,
Georg-Au\-gust-Uni\-versi\-t{\"a}t G{\"o}t\-ting\-en,
Bunsenstra{\ss}e 3-5,
D-37073 G{\"o}ttingen, 
Germany.}
\email{kyed@uni-math.gwdg.de}
\urladdr{www.uni-math.gwdg.de/kyed}

\thanks{Research supported by The Danish Council for Independent Research $\mid$ Natural Sciences}

\begin{document}
\begin{abstract}
We prove a Delorme-Guichardet type theorem for discrete quantum groups expressing property (T) of the quantum group in question in terms of its first cohomology groups. As an application, we show that the first $L^2$-Betti number of a discrete property (T) quantum group vanishes.

\end{abstract}
\maketitle
\section{Introduction}
The notion of property $\T$ was introduced by Kazhdan in his influential paper \cite{kazhdan} and has since then played a prominent role in a variety of mathematical disciplines, including  topology, ergodic theory and operator algebras.  Over the years the definition has been generalized to different operator algebraic settings, for instance  by Connes and Jones for the class of $\I\I_1$-factors  in \cite{connes-jones} and by Bekka for tracial $C^*$-algebras in \cite{bekka-T-for-C}. Recently Fima introduced property $\T$ in the context of discrete quantum groups; a class of operator algebras not necessarily arising from groups, but still carrying some of the extra structure present in group $C^*$-algebras or group von Neumann algebras. The present paper is devoted to the study of this notion of property $\T$. Before stating our main results, we set the stage by briefly discussing the definitions and a few classical results concerning property $\T$ for discrete groups; for details the reader is referred to monograph \cite{BHV} by Bekka, de la Harpe and Valette.  Recall that a discrete, countable group $\Gamma$ has property $\T$  if any unitary representation of $\Gamma$ with almost invariant vectors has a non-zero invariant vector.   
That is, if a unitary representation $\pi\colon \Gamma \to B(H)$ admits a sequence of unit vectors $\xi_n\in H$ such that $\|\pi(\gamma)\xi_n-\xi_n  \|\to 0$ for every $\gamma\in \Gamma$, then there exists a non-zero vector $\xi\in H$ such that $\pi(\gamma)\xi=\xi$ for every $\gamma\in \Gamma$.
One reason why property $\T$ is an important notion is that it allows many different descriptions. Firstly, it can be described using the positive definite functions functions on $\Gamma$ by means of the following theorem.

\begin{thm*}[Akemann-Walter, \cite{akemann-walter}]
The group $\Gamma$ has property $\T$ if and only if any sequence of  positive definite functions $\varphi_n\colon \Gamma\to \CC$  converging pointwise to 1 and with $\varphi_n(e)=1$, converges uniformly to the constant function 1. 
\end{thm*}

Property $\T$ can also be described in terms of the first cohomology of $\Gamma$ which, among other things, provides a link between property $\T$ and Serre's property (FA).  The precise cohomological description is given by the celebrated Delorme-Guichardet theorem.

\begin{thm*}[Delorme-Guichardet, \cite{delorme,guichardet}]	
The group $\Gamma$ has property $\T$ if and only if the first group cohomology $H^1(\Gamma,H)$ vanishes for all Hilbert spaces $H$ carrying a unitary $\Gamma$-action.
\end{thm*}
Both of these results have analogues in the $\I\I_1$-factor setting  \cite{peterson-T}, but statements and proofs are considerably more involved than in the case of discrete groups. The main purpose of the present paper is to show how the classical results mentioned above can be generalized to the quantum group context in a way that is somewhat closer to the original results than the generalizations obtained in the general setting of von Neumann algebras with property $\T$. If $\hat{\GG}$ is a discrete quantum group and $\GG$ is its compact dual, we denote by $(\Pol(\GG),\Delta, S,\varps)$ the associated Hopf $*$-algebra of matrix coefficients and by $C(\GG_\max)$ the universal $C^*$-completion of $\Pol(\GG)$. These objects will be introduced and discussed in greater detail in Section \ref{qgrp-preliminaries} where we also elaborate on Fima's definition of property $\T$ and the results obtained in \cite{fima-prop-T}. Our first main result (Theorem \ref{convergence-thm})  is an analogue for quantum groups of the result of Akemann and Walter.

\begin{thm*}
The discrete quantum group $\hat{\GG}$ has property $\T$ if and only if any net of states  $\varphi_i\colon C(\GG_\max)\to\CC$ converging pointwise to the counit $\varps$ converges in the uniform norm.
\end{thm*}
Secondly, we prove in Theorem \ref{delorme-guichardet} the following quantum group version of the Delorme-Guichardet theorem.
\begin{thm*}

The discrete quantum group $\hat{\GG}$ has property $\T$ if and only if the following holds: for every $*$-representation $\pi\colon\Pol(\GG)\to B(H)$ on a Hilbert space $H$ the first Hochschild cohomology of $\Pol(\GG)$ with values in the  bimodule $_\pi H_\varps$ vanishes.
\end{thm*}
The relevant definitions concerning the first Hochschild cohomology will be given in Section \ref{cocycle-section}. Along the way we also obtain (see Theorem \ref{delorme-guichardet}) a characterization of property $\T$ in terms of conditionally negative functionals $\psi\colon\Pol(\GG)\to \CC$  that parallels the classical description stating that a discrete group $\Gamma$ has property $\T$ if and only if every conditionally negative definite function $\psi\colon \Gamma\to \RR$ is bounded. Finally, as an application we obtain  in Corollary \ref{beta-et=0} the following generalization of a well known result (see \cite{bekka-valette-group-cohomology}) for groups.
\begin{cor*}
If $\hat{\GG}$ has property $\T$ then its first $L^2$-Betti number vanishes.
\end{cor*}
The paper is organized as follows.
\vspace{0.3cm}
\paragraph{\emph{Structure.}}
The first section provides the reader with the necessary background concerning the theory of compact quantum groups, their discrete duals and  the definition of property $\T$ for discrete quantum groups.  In Section \ref{dual-point-section} we show how property $\T$ of a discrete quantum group can be described in terms of its dual compact quantum group and use this description to give a spectral interpretation of property $\T$. Section \ref{pos-def-section} is devoted to the proof of the characterization of property $\T$ in terms of  states on the associated universal $C^*$-algebra, and in Section \ref{cocycle-section} the proof of the Delorme-Guichardet theorem is given. In the sixth and final section we show how the results obtained can be used to derive information about the $L^2$-invariants of the quantum group in question.

\vspace{0.3cm}
\paragraph{\emph{Notation.}}
Throughout the paper, the symbol $\odot$ will be used to denote algebraic tensor products while the symbol $\htens$ will be used to denote tensor products of Hilbert spaces and von Neumann algebras. All tensor products between $C^*$-algebras are assumed minimal/spatial and these will be denoted by the symbol $\tens$. Hilbert spaces are assumed to be complex and their inner products to be linear in the first variable. Furthermore, $*$-representations of unital algebras on Hilbert spaces are implicitly assumed to be unit-preserving.

\vspace{0.3cm}
\paragraph{\emph{Acknowledgements.}}
The author wishes to thank Pierre Fima, Ryszard Nest, Jesse Peterson and Roland Vergnioux for discussions revolving around the notion of property $\T$. The work presented was initiated during the authors stay at the Hausdorff Research Institute for Mathematics  and it is a pleasure to thank the organizers, Wolfgang L{\"u}ck and Nicolas Monod, as well as the participants of the trimester program on Rigidity.

\section{Preliminaries on quantum groups}\label{qgrp-preliminaries}

We choose here the approach to compact quantum groups developed by Woro\-no\-wicz in \cite{woronowicz-pseudo}, \cite{woronowicz} and \cite{wor-cp-qgrps}. Thus, a compact quantum group $\GG$ consists of a (not necessarily commutative) separable, unital $C^*$-algebra $C(\GG)$ together with a unital, coassociative $*$-homomorphism $\Delta\colon C(\GG)\to C(\GG)\tens C(\GG)$ satisfying a certain density condition. The map $\Delta$ is referred to as the comultiplication. Such a quantum group possesses a unique Haar state; i.e.~a state $h\colon C(\GG)\to \CC$ such that
\[
(\id\tens h)\Delta a=h(a)1=(h\tens \id)\Delta (a) \text{ for every } a\in C(\GG).
\]
The GNS-construction applied to the Haar state yields a separable Hilbert space $L^2(\GG)$ together with a $*$-representation $\lambda\colon C(\GG)\to B(L^2(\GG))$ and a linear map $\Lambda\colon C(\GG)\to L^2(\GG)$ with dense image. In general, the Haar state need not be faithful and hence the \emph{left regular representation} $\lambda$ might have a kernel and we denote by $C(\GG_\red)$ the image $\lambda(C(\GG))$. This $C^*$-algebra inherits a quantum group structure  from $\GG$ and the comultiplication $\Delta_\red$ on  $C(\GG_\red)$ is implemented by the so-called multiplicative unitary $W\in B(L^2(\GG)\htens L^2(\GG))$ given by
\[
W^*(\Lambda(a)\tens\Lambda(b))=\Lambda\tens\Lambda(\Delta(b)(a\tens 1)).
\]
The statement that $W$ implements $\Delta_\red$ means that $\Delta_\red(\lambda(a))=W^*(1\tens \lambda(a))W$ for every $a\in C(\GG)$. One notes that the right hand side of this formula also makes sense if $\lambda(a)$ is replaced with any $T\in B(L^2(\GG))$, and it turns out that the enveloping von Neumann algebra $L^\infty(\GG)=C(\GG_\red)''$ is turned into a compact von Neumann algebraic quantum group (see \cite{kustermans-vaes}) when endowed with this map as comultiplication. \\

A unitary corepresentation of $\GG$ on a Hilbert space $H$ is a unitary element $u\in \mathcal{M}(\K(H)\tens C(\GG))$ such that $(\id\tens \Delta)u=u_{(12)}u_{(13)}$. Here $\K(H)$ denotes the compact operators on $H$, $\mathcal{M}(\K(H)\tens C(\GG))$ is the multiplier algebra of $\K(H)\tens C(\GG)$ and the subscripts $(12)$ and $(13)$ are the standard leg-numbering notation.  Representation theoretic notions from the theory of compact groups, such as direct sums, tensor products, intertwiners and irreducibility, have natural counterparts in the corepresentation theory for compact quantum groups. In particular the following important theorem holds true.

\begin{thm}[Woronowicz]
Any irreducible unitary corepresentation of $\GG$ is finite dimensional and an arbitrary unitary corepresentation decomposes as a direct sum of irreducible ones.
\end{thm}
We denote by $\Irred(\GG)$ the set of equivalence classes of irreducible, unitary corepresentations of $\GG$. The separability assumption on $C(\GG)$ together with the quantum Peter-Weyl theorem \cite{wor-cp-qgrps} ensures that $\Irred(\GG)$ is a countable set. We label its elements by an auxiliary  countable set $I$ and choose for each $\alpha\in I$ a Hilbert space $H^\alpha$ and a concrete representative $u^\alpha\in B(H^\alpha)\tens C(\GG)$. Abusing notation slightly, we shall often identify the index $\alpha$ with the corresponding class of $u^\alpha$. Fix an orthonormal basis $\{e_1,\dots, e_{n_\alpha}\}$ for $H^\alpha$ and consider the corresponding functionals $\omega_{ij}\colon B(H^\alpha)\to \CC$ given by $\omega_{ij}(T)=\ip{Te_i}{e_j}$. The matrix coefficients of $u^\alpha$, relative to the chosen basis, are then defined as
\[
u_{ij}^\alpha=(\omega_{ji}\tens \id)u^\alpha \in C(\GG).
\]
It turns out that these matrix coefficients are linearly independent and that their linear span constitutes a dense $*$-subalgebra $\Pol(\GG)$ of $C(\GG)$. Furthermore, the comultiplication descends to a comultiplication $\Delta\colon \Pol(\GG)\to \Pol(\GG)\odot \Pol(\GG)$ and with this comultiplication $\Pol(\GG)$ becomes a Hopf $*$-algebra; i.e.~there exists an antipode $S\colon \Pol(\GG)\to \Pol(\GG)$ as well as a counit $\varps\colon \Pol(\GG)\to \CC$ satisfying the usual Hopf $*$-algebra relations \cite{klimyk}. The fact that $\Pol(\GG)$ is spanned by matrix coefficients arising from finite dimensional, unitary corepresentations also ensures that the relation
\[
\|a\|_\max=\sup\{\|\pi(a)\|\mid \pi\colon \Pol(\GG)\to B(H) \text{ a cyclic $*$-representation}\}
\]
defines a $C^*$-norm $\|\cdot\|_\max$ on $\Pol(\GG)$ which dominates any other $C^*$-norm. The $C^*$-completion of $\Pol(\GG)$ with respect to this norm is called the universal $C^*$-algebra associated with $\GG$ and is denoted $C(\GG_\max)$. By definition of $\|\cdot\|_\max$, the comultiplication extends to a comultiplication $\Delta_\max\colon C(\GG_\max)\to C(\GG_\max)\tens C(\GG_\max)$ turning $C(\GG_\max)$ into a compact quantum group. Note that the $*$-representations of $C(\GG_\max)$ are in one-to-one correspondence with the $*$-representations of $\Pol(\GG)$ via restriction/extension.

\begin{ex}
The fundamental example of a compact quantum group, on which the general definition is modeled, is obtained by considering a compact, second countable, Hausdorff topological group $G$ and its commutative $C^*$-algebra $C(G)$ of continuous, complex valued functions. In this case the comultiplication is the Gelfand dual of the multiplication map $G\times G\to G$ and the Haar state is given by integration against the unique Haar probability measure $\mu$ on $G$. The GNS-space therefore identifies with $L^2(G,\mu)$ and the $*$-representation $\lambda$ with the action of $C(G)$ on $L^2(G,\mu)$ by pointwise multiplication. Similarly, the von Neumann algebra identifies with $L^\infty(G,\mu)$ and the Hopf $*$-algebra becomes the subalgebra of $C(G)$ generated by matrix coefficients arising from irreducible, unitary representations of $G$. The antipode is the Gelfand dual of the inversion map and the counit is given by evaluation at the neutral element in $G$.
\end{ex}
In the previous example there is no real difference between the reduced and universal version of the compact quantum group. The next example, however, will illustrate this difference more clearly.

\begin{ex}
Consider a countable, discrete group $\Gamma$. Denote by $C^*_\red(\Gamma)$ its reduced group $C^*$-algebra acting on $\ell^2(\Gamma)$ via the left regular representation and define a comultiplication on group elements by $\Delta_\red\gamma=\gamma\tens\gamma$. This turns $C^*_\red(\Gamma)$ into a compact quantum group whose Haar state is given by the natural trace on $C^*_\red(\Gamma)$. Hence the GNS-space and GNS-representation can be identified, respectively, with $\ell^2(\Gamma)$ and the left regular representation,  and the enveloping von Neumann algebra is therefore nothing but the group von Neumann algebra $\L(\Gamma)$. Each element in $\Gamma$ is a one-dimensional corepresentation for this quantum group and the Hopf $*$-algebra therefore identifies with the complex group algebra $\CC\Gamma$. Thus, the universal $C^*$-algebra is, by definition, equal to the maximal group $C^*$-algebra $C^*_\max(\Gamma)$.
\end{ex}

\begin{rem}
The three $C^*$-algebras $C(\GG), C(\GG_\red)$ and $C(\GG_\max)$, together with their comultiplications, can be thought of as ``different pictures of the same quantum group'', each having its advantages and disadvantages. For instance, whereas the Haar state is always faithful on $C(\GG_\red)$ this is in general not the case on $C(\GG_\max)$ and, conversely, the counit is always well defined on all of $C(\GG_\max)$ but not necessarily on $C(\GG_\red)$. The latter difference is the fundamental observation leading to the notion of (co-)amenability for quantum groups as studied in \cite{murphy-tuset}.
\end{rem}

Any compact quantum group $\GG$ has a dual quantum group $\hat{\GG}$ of so-called discrete type. As in the compact case, $\hat{\GG}$ comes with both a $C^*$-algebra $c_0(\hat{\GG})$ and a von Neumann algebra $\ell^\infty(\hat{\GG})$ defined, respectively, as
\[
c_0(\hat{\GG})=\bigoplus_{\alpha\in I}^{c_0}B(H^\alpha)\ \text{ and } \ \ell^\infty(\hat{\GG})=\prod_{\alpha\in I}^{\ell^\infty}B(H^\alpha).
\]
In the discrete picture we will primarily be working with the von Neumann algebra $\ell^\infty(\hat{\GG})$, which is endowed with a natural comultiplication $\hat{\Delta}\colon \ell^\infty(\hat{\GG})\to \ell^\infty(\hat{\GG})\htens \ell^\infty(\hat{\GG})$ arising from the quantum group structure on $\GG$. 
Since $c_0(\hat{\GG})$ is a direct sum of finite dimensional $C^*$-algebras we have isomorphisms
\[
\ell^{\infty}(\hat{\GG})\htens B(H)\simeq \mathcal{M}(c_0(\hat{\GG})\tens B(H))\simeq \prod_{\alpha\in I}^{\ell^\infty}B(H^\alpha)\tens B(H)
\]
for any Hilbert space $H$. For an element $T\in \ell^{\infty}(\hat{\GG})\htens B(H) $ we will denote by $(T^\alpha)_{\alpha\in I}$ the corresponding element in $\prod_{\alpha\in I}^{\ell^\infty}B(H^\alpha)\tens B(H)$ and in the sequel we will freely identify $T$ and $(T^\alpha)_{\alpha\in I}$. By a unitary corepresentation of $\hat{\GG}$ on a Hilbert space $H$ we shall mean a unitary operator $V\in \ell^\infty(\hat{\GG})\htens B(H)$ satisfying
\[
(\hat{\Delta}\tens \id)V=V_{(23)}V_{(13)}.
\]
Consider the so-called \emph{universal bicharacter} $V_\max = (u^\alpha)_{\alpha\in I}\in \prod_{\alpha\in I}^{\ell^\infty}B(H_\alpha)\tens C(\GG_\max)$. This unitary encodes the duality between $\GG$ and $\hat{\GG}$ in the following sense: for every unitary corepresentation $V\in \ell^\infty(\hat{\GG})\htens B(H)$ of $\hat{\GG}$ there exists 
a unique $*$-representation $\pi_V\colon C(\GG_\max)\to B(H)$ such that
\[
(\id\tens \pi_V)u^\alpha=V^\alpha \ \text{ for each } \alpha\in I.
\]
Conversely, every $*$-representation $\pi\colon C(\GG_\max)\to B(H)$ defines a unitary corepresentation of $\hat{\GG}$ on $H$ by the above relation. 
See \cite{soltan-woronowicz} for details.\\

	As mentioned in the introduction, the notion of property $\T$ was recently introduced in the quantum group setting by Fima \cite{fima-prop-T} and the definition is as follows.

\begin{defi}[Fima]
Let $\hat{\GG}$ be a discrete quantum group and consider a unitary corepresentation $V=(V^\alpha)_{\alpha\in I}\in \ell^\infty(\hat{\GG})\htens B(H)$ of $\hat{\GG}$ on a Hilbert space $H$. 
\begin{itemize}
\item[(i)] A vector $\xi\in H$ is said to be invariant if $V^\alpha(\eta\tens\xi)=\eta\tens\xi$ for all $\alpha\in I$ and $\eta\in H^\alpha$. 
\item[(ii)] For a finite, non-empty subset $E\subseteq \Irred(\GG)$ and a $\delta>0$ a non-zero vector $\xi\in H$ is called $(E,\delta)$-invariant if
$
\|V^\alpha(\eta\tens\xi)-\eta\tens \xi \|<\delta\|\eta\|\|\xi\|$  for all  $\alpha\in E$ and all  $\eta\in H^\alpha$,
and $V$ is said to have almost-invariant vectors if it has an $(E,\delta)$-invariant vector for each finite, non-empty $E\subseteq \Irred(\GG)$ and each $\delta>0$.
\item[(iii)] The discrete quantum group $\hat{\GG}$ is said to have property $\T$ if any unitary corepresentation of $\hat{\GG}$ with almost invariant vectors has a non-zero invariant vector.
\end{itemize}
\end{defi}
\begin{rem}
For notational smoothness we will adopt the convention that, unless explicitly specified otherwise, subsets $E$ of $\Irred(\GG)$ are always both finite and non-empty.
\end{rem}
\begin{rem}
The study of property $\T$ for quantum groups began before the paper \cite{fima-prop-T}. In \cite{petrescu} property $\T$ was studied in the setting of Kac algebras and in \cite{bedos-conti-tuset} it was introduced for the the class of algebraic quantum groups. As is shown in \cite{exotic-norms}, these different notions all agree with Fima's definition in the case of a discrete quantum group.

\end{rem}

The main results in \cite{fima-prop-T} are summarized in the following theorem.
\begin{thm}[Fima]\label{fima-thm}
If $\hat{\GG}$ is a discrete quantum group with property $\T$ then the following holds.
\begin{itemize}
\item[(i)] The quantum group is automatically of Kac type; i.e.~the Haar state $h\colon C(\GG)\to \CC$ is a trace.
\item[(ii)] The discrete quantum group is finitely generated; i.e.~the corepresentation category $\Corep(\GG)$ of the compact dual is a finitely generated tensor category.
\item[(iii)] The quantum group allows Kazhdan pairs; i.e.~for every finite subset $E\subseteq \Irred(\GG)$ generating the corepresentation category and containing the trivial corepresentation there exists a $\delta>0$ such that whenever $V$ is a unitary corepresentation of $\hat{\GG}$ having an $(E,\delta)$-invariant vector, then $V$ has a non-trivial invariant vector.
\end{itemize}
Moreover, if $\hat{\GG}$ is an infinite, discrete quantum group such that $L^\infty(\GG)$ is a factor then $\hat{\GG}$ has property $\T$ iff $L^\infty(\GG)$ is a type $\I\I_1$-factor with property $\T$ in the sense of Connes-Jones \cite{connes-jones}.

\end{thm}
\begin{rem}
Concrete non-cocommutative examples of quantum groups with property $\T$ was constructed in \cite[Example 3.1]{fima-prop-T} by twisting the comultiplication on $\widehat{SL_{n}(\ZZ)}$  by a 2-cocycle. Using \cite[Proposition 6.1]{exotic-norms} it is not difficult to see that one can also obtain examples, of an admittedly somewhat trivial nature, by considering quantum groups of the form $\widehat{\GG \times  \HH}$ where $\hat{\GG}$ is a discrete (quantum) group with property $\T$ and $\HH$ is any finite quantum group.  
\end{rem}

\section{Property (T) from the dual point of view}\label{dual-point-section}
In this section we reformulate property $\T$ for discrete quantum groups in terms of their compact duals and use this description to give a spectral characterization of property $\T$.  In the compact setting it is natural to consider the following notions of invariance and almost invariance.
\begin{defi}
Let $\pi\colon \Pol(\GG)\to B(H) $ be a $*$-representation. A vector $\xi\in H$ is said to be invariant if  $\pi(a)\xi=\varps(a)\xi$ for all $a\in \Pol(\GG)$. If a non-zero invariant vector exists then $\pi$ is said to contain the counit. For a  subset $E\subseteq \Irred(\GG)$ and $\delta >0$ a vector $\xi\in H$ is said to be $(E,\delta)$-invariant if 
\[
\|\pi(u_{ij}^\alpha)\xi-\varps(u_{ij}^\alpha)\xi \|<\delta\|\xi\|,
\]
for all $\alpha\in E$ and $i,j\in\{1,\dots, n_\alpha\}$. The $*$-representation $\pi$ is said to have almost invariant vectors if it allows a non-zero $(E,\delta)$-invariant vector for every finite  $E\subseteq \Irred(\GG)$ and every $\delta>0$.
\end{defi}

\begin{rem}
Since the set $\{u_{ij}^\alpha\mid \alpha\in I, 1\leq i,j \leq n_\alpha\}$ spans $\Pol(\GG)$ linearly it is not difficult to see that a $*$-representation $\pi\colon \Pol(\GG)\to B(H)$ has almost invariant vectors iff there exists a sequence of unit vectors $\xi_n\in H$ such that 
\[
\lim_{n\to\infty} \|\pi(a)\xi_n-\varps(a)\xi_n\|=0 
\]
for every $a\in \Pol(\GG)$.
\end{rem}

The following proposition contains the translation of property $\T$ from the discrete to the compact picture.

\begin{prop}\label{dual-invariance}
Let $\hat{\GG}$ be a discrete quantum group and consider a unitary corepresentation $V\in \ell^\infty(\hat{\GG})\htens B(H)$ as well as the the corresponding $*$-rep\-res\-entation $\pi_V\colon \Pol(\GG)\to B(H)$. 
Let furthermore $E\subseteq \Irred(\GG)$ and $\delta>0$ be given and define  $K_E=\maxx\{n_\alpha \mid \alpha\in E\}$.
Then the following holds.
\begin{itemize}
\item[(i)] A vector $\xi\in H$ is $V$-invariant if and only if it is $\pi_V$-invariant.
\item[(ii)] If  $\xi\in H$ is $(E,\delta)$-invariant for $V$ then it is also $(E, \delta)$-invariant for $\pi_V$.
\item[(iii)] If $\xi\in H$ is $(E,\delta)$-invariant for $\pi_V$ then it is $(E,K_E\delta)$-invariant for $V$.
\end{itemize}
Thus, $V$ has almost-invariant vectors if and only if $\pi_V$ has almost invariant vectors. 
In particular, the discrete quantum group $\hat{\GG}$ has property $\T$ iff  any $*$-representation $\pi\colon \Pol(\GG)\to B(H)$ with almost invariant vectors has a non-zero invariant vector.
\end{prop}

\begin{proof}
Consider the fixed basis $\{e_1,\dots,e_{n_\alpha}\}$ of $H_\alpha$ and the corresponding functionals $e_i'\colon H_\alpha\to\CC$ given by $e_i'(x)=\ip{x}{e_i}$. A vector $\xi\in H$ is $V$-invariant exactly when $V^\alpha(e_j\tens \xi)=e_j\tens\xi$ for any $\alpha\in I$ and $j\in \{1,\dots, n_\alpha\}$. This in turn holds iff
\[
(e_i'\tens\id)V^\alpha(e_j\tens \xi)= e_i'(e_j) \xi \ \text{ for all } \alpha\in I  \text{ and } i,j\in \{1,\dots, n_\alpha\}.
\]
Keeping in mind that $\varps(u_{ij}^\alpha)=\delta_{ij}$, the above equation translates into
\[
\pi_V(u_{ij}^\alpha)\xi= \varps(u_{ij}^\alpha)\xi \ \text{ for all } \alpha\in I  \text{ and } i,j\in \{1,\dots, n_\alpha\},
\]
which is equivalent to $\xi$ being $\pi_V$-invariant since the matrix coefficients constitute a linear basis for $\Pol(\GG)$. This proves (i). To prove (ii), fix $E\subseteq \Irred(\GG)$ and $\delta>0$ and assume that $\xi\in H$ is an $(E,\delta)$-invariant unit vector for $V$. Then for each $\alpha\in E$ and $i,j \in \{1,\dots,n_\alpha\}$ we have 
\[
\|\pi_V(u_{ij}^\alpha)\xi-\varps(u_{ij}^\alpha)\xi\|=\|(e_i'\tens\id)(V^\alpha(e_j\tens \xi)-e_j\tens\xi)\|\leq \|V^\alpha(e_j\tens \xi)-e_j\tens\xi\|<\delta,
\]  
as desired. To prove (iii), assume that $\xi\in H$ is an $(E,\delta)$-invariant unit vector for $\pi_V$. For each $\alpha\in E$ and  $j\in \{1,\dots, n_\alpha\}$ we then get
\begin{align*}
\|V^\alpha(e_j\tens\xi)-e_j\tens\xi\|^2 &=\sum_{i=1}^{n_\alpha}\|(e_i'\tens\id)\left ( V^\alpha(e_j\tens\xi)-e_j\tens\xi\right )\|^2\\
&=\sum_{i=1}^{n_\alpha} \|\pi_V(u_{ij}^\alpha)\xi-\varps(u_{ij}^\alpha)\xi\|^2<n_\alpha\delta^2. 
\end{align*}  
Hence for $\eta=\sum_{i=1}^{n_\alpha}\eta_ie_i\in H_\alpha$ we get by H\"older's inequality
\begin{align*}
\|V^\alpha(\eta\tens \xi)-\eta\tens\xi\| &\leq \sum_{i=1}^{n_\alpha}|\eta_i|\|V^\alpha(e_i\tens\xi)-e_i\tens\xi\|\\
&< \|\eta\|_1\sqrt{n_\alpha}\delta\\
&\leq \|\eta\|_2 n_\alpha\delta, 
\end{align*}
which shows that $\xi$ is $(E, K_E\delta)$-invariant for $V$.
\end{proof}

Similarly, the existence of Kazhdan pairs also translates to the dual picture.
\begin{cor}\label{kazhdan-pair-cor}
Let $\hat{\GG}$ have property $\T$ and let $E\subseteq \Irred(\GG)$ be a finite subset containing the trivial corepresentation $1$ which generates the corepresentation category of $\GG$. Then there exists $\delta>0$ such that any $*$-representation $\pi\colon\Pol(\GG)\to B(H)$ having an $(E,\delta)$-invariant vector has a non-zero invariant vector.
\end{cor}
\begin{proof}
Let a $*$-representation $\pi\colon \Pol(\GG)\to B(H)$ be given. Denote by $V$ the corresponding corepresentation of $\hat{\GG}$ on $H$ and choose $\delta>0$ such that $(E,\delta)$ is a Kazhdan pair for $\hat{\GG}$. If we put $K_E=\maxx\{n_\alpha\mid \alpha\in E\}$ and  $\xi\in H$ is an $(E,K_E^{-1}\delta)$-invariant vector for $\pi$ then, by Proposition \ref{dual-invariance} (iii), $\xi$ is an $(E,\delta)$-invariant vector for $V$. Hence $V$ allows a non-zero invariant vector which is then also invariant for $\pi$ by Proposition \ref{dual-invariance} (i). 
\end{proof}

\begin{rem}
In the sequel we will primarily work with the $*$-representations of $\Pol(\GG)$ instead of the corepresentations of $\hat{\GG}$ and thus a \emph{Kazhdan pair} for $\hat{\GG}$ is going to mean a pair as described in Corollary \ref{kazhdan-pair-cor}.
\end{rem}

Recall from Theorem \ref{fima-thm} that a discrete property $\T$ quantum group is automatically finitely generated. Consider now any finitely generated, discrete quantum group $\hat{\GG}$ and let $E\subseteq \Irred(\GG)$ be a finite generating set for $\Corep(\GG)$ containing the trivial corepresentation. For each $\alpha\in E$ and $i,j\in\{1,\dots,n_\alpha\}$  define $x_{ij}^\alpha=u_{ij}^\alpha-\varps(u_{ij}^\alpha)1$ and put $X_E=\sum_{\alpha\in E, i,j}x_{ij}^{\alpha*}x_{ij}^\alpha$.  Property $\T$ can then be read of the element $X_E$ by means of the following result.

\begin{thm}\label{spectral-T}
The discrete quantum group $\hat{\GG}$ has property $\T$ if and only if zero is not in the spectrum of $\pi(X_E)$ for any $*$-representation $\pi\colon \Pol(\GG)\to B(H)$ not containing the counit.
\end{thm}

\begin{proof}
Assume that $\hat{\GG}$ has property $\T$ and let $\pi\colon \Pol(\GG)\to B(H)$ be a $*$-representation such that $\pi(X_E)$ is not bounded away from zero. Then there exists a sequence $(\xi_k)_{k\in \NN}$ in the unit ball of $H$ such that $\pi(X_E)\xi_k\to 0$. Hence
\begin{align*}
0&=\lim_k\ip{\pi(X_E)\xi_k}{\xi_k}\\
&=\lim_{k\to\infty}\sum_{\alpha\in E}\sum_{i,j=1}^{n_\alpha}\ip{\pi(x_{ij}^\alpha)^*\pi(x_{ij}^\alpha)\xi_k}{\xi_k}\\
&=\lim_{k\to\infty}\sum_{\alpha\in E}\sum_{i,j=1}^{n_\alpha}\|\pi(u_{ij}^\alpha)\xi_k-\varps(u_{ij}^\alpha)\xi_k \|^2.
\end{align*}
For a suitable $\delta>0$ the pair $(E,\delta)$ is a Kazhdan pair for $\hat{\GG}$ and therefore the above convergence forces $\pi$ to have a non-trivial invariant vector; hence $\pi$ contains $\varps$.  Conversely, assume that $\hat{\GG}$ does not have property $\T$. Then there exists a $*$-representation $\pi\colon \Pol(\GG)\to B(H)$ with almost invariant vectors, but without non-zero invariant vectors. In particular we may find a sequence of unit vectors $(\xi_k)_{k\in \NN}$ in $H$ such that 
\[
\lim_{k\to \infty}\sum_{\alpha\in E}\sum_{i,j=1}^{n_\alpha}\|\pi(u_{ij}^\alpha)\xi_k-\varps(u_{ij}^\alpha)\xi_k\|^2=0.
\]
On the other hand we have
\[
\sum_{\alpha\in E}\sum_{i,j=1}^{n_\alpha}\|\pi(u_{ij}^\alpha)\xi_k-\varps(u_{ij}^\alpha)\xi_k\|^2=\ip{\pi(X_E)\xi_k}{\xi_k}=\|\pi(X_E)^{\frac12}\xi_k\|^2,
\]
and hence zero is in the spectrum of $\pi(X_E)^{\frac12}$. Thus $\pi$ is a $*$-representation not containing $\varps$ such that $\pi(X_E)^{\frac12}$, and hence also $\pi(X_E)$, is not invertible.

\end{proof}

\begin{rem}
The above spectral characterization of property $\T$ should be compared with the Kesten condition for coamenability \cite{banica-subfactor}, which states that $\GG$ is coamenable iff zero is in the spectrum of $\lambda(X_E)$. Theorem \ref{spectral-T} is an extension of a result for groups due to de la Harpe, Robertson and Valette \cite{valette-spectrum-of-sum}.
\end{rem}

As in the classical situation, we also get a version of property $\T$ with ``continuity constants''. 
\begin{prop}\label{continuity-prop}
A discrete quantum group $\hat{\GG}$ has property $\T$ if and only if one of  the following two condition holds.
\begin{itemize}
\item[(i)] For every $\delta>0$ there exists $E_0\subseteq \Irred(\GG)$ and $\delta_0>0$ such that any $*$-representation $\pi\colon\Pol(\GG)\to B(H)$ with an $(E_0,\delta_0\delta)$-invariant vector $\xi\in H$ has an invariant vector $\eta\in H$ such that $\|\xi-\eta\|<\delta \|\xi\|$. 
\item[(ii)] For every $\delta>0$ there exist $E_0\subseteq \Irred(\GG)$ and $\delta_0>0$ such that any $*$-representation $\pi\colon\Pol(\GG)\to B(H)$ with an $(E_0,\delta_0)$-invariant unit vector $\xi\in H$ has an invariant vector $\eta\in H$ such that $\|\xi-\eta\|<\delta$.
\end{itemize}
\end{prop}
The proof of the proposition is basically identical to the corresponding proof in the group case \cite[Proposition 1.1.19]{BHV}, but we include the short argument for the sake of completeness.
\begin{proof}
Assume that $\hat{\GG}$ has property $\T$ and let $\delta>0$ be given. Choose a Kazhdan pair $(E_0,\delta_0)$ for $\hat{\GG}$ and consider a $*$-representation $\pi\colon\Pol(\GG)\to B(H)$. Denote by $P\in B(H)$ the projection onto the closed subspace of invariant vectors. Assume furthermore that $\xi$ is an $(E_0,\delta_0\delta)$-invariant vector and decompose $\xi$ as $\xi=\xi'+\xi''$ with $\xi'=P\xi$ and $\xi''=(1-P)\xi$. Since $P(H)^\perp$ does not have non-zero invariant vectors and $(E_0,\delta_0)$ is a Kazhdan pair there must exist $\beta\in E_0$ and $k,l\in \{1,\dots, n_\beta\}$ such that
\[
\|\pi(u_{kl}^\beta)\xi''-\varps(u_{kl}^\beta)\xi''\|\geq \delta_0 \|\xi''\|.
\]
Using that $\xi$ is $(E_0,\delta_0\delta)$-invariant we get
\[
\delta_0\delta\|\xi\| > \|\pi(u_{kl}^\beta)\xi-\varps(u_{kl}^{\beta})\xi\|= \|\pi(u_{kl}^\beta)\xi''-\varps(u_{kl}^{\beta})\xi''\|\geq \delta_0\|\xi''\|,
\]
and hence that $\delta\|\xi\|> \|\xi''\|$. Putting $\eta=\xi'$ we get 
\[
\|\xi-\eta\|=\|\xi-\xi'\|=\|\xi''\|< \delta\|\xi\|
\]
and (i) follows. To prove (ii),  let $\delta>0$ be given and assume without loss of generality that  $\delta\leq 1$. Choose a Kazhdan pair $(E_0,\delta_0')$ for $\hat{\GG}$ and put $\delta_0=\delta_0'\delta$.  Then $(E_0,\delta_0)$ is also a Kazhdan pair and from the proof of (i) we have that the pair $(E_0,\delta_0)$ satisfies the claim. That (i) and (ii) both imply property $\T$ is clear.

\end{proof}

\section{\texorpdfstring{Property (T) in terms of states on the universal $C^*$-algebra}{Property (T) in terms of states on the universal C*-algebra}}\label{pos-def-section}

As already mentioned in the introduction, a discrete group $\Gamma$ has property $\T$ exactly when every sequence of normalized, positive definite functions on $\Gamma$ converging pointwise to $1$ actually converges uniformly to 1.
Recall that a function $\varphi\colon \Gamma\to \CC$ is called normalized if $\varphi(e)=1$ and positive definite if
\[
\sum_{i,j=1}^n \bar{\alpha}_i\alpha_j \varphi(\gamma_i^{-1}\gamma_j)\geq 0 \ \text{ for all } n\in \NN, \gamma_1,\dots,\gamma_n\in \Gamma \text{ and } \alpha_1,\dots, \alpha_n\in \CC.
\]
Hence, there is a one-to-one correspondence between normalized, positive definite functions on $\Gamma$  and states on the universal group $C^*$-algebra $C^*_\max(\Gamma)$.  Having this correspondence in mind, the following theorem generalizes the classical result.

\begin{thm}\label{convergence-thm}
A discrete quantum group $\hat{\GG}$ has property $\T$ if and only if any net of states on $C(\GG_\max)$ converging pointwise to the counit $\varps$ converges in the uniform norm.
\end{thm}
Here \emph{the uniform norm} is the norm on the state space of $C(\GG_\max)$ given by
\[
\|\varphi\|=\sup\{|\varphi(a)|\mid \|a\|_\max\leq 1\},
\]
and convergence in this norm will be referred to as uniform convergence. The fact that property $\T$ implies the convergence property was proved independently by Fima (private communication) in the dual picture.  For the proof of Theorem \ref{convergence-thm} we need the following lemma.

\begin{lem}\label{unitary-lem}
Let $\hat{\GG}$ be a discrete quantum group and let $\delta>0$ and a $*$-repre\-sen\-tation $\pi\colon C(\GG_\max)\to B(H)$ be given. If $\xi\in H$ is a unit vector such that $\|\pi(v)\xi-\varps(v)\xi\|\leq \delta$ for every unitary $v\in C(\GG_\max)$ then there exists an invariant vector $\eta\in H$ such that $\|\xi-\eta\|\leq \delta$.
\end{lem}
The proof of Lemma \ref{unitary-lem} is inspired by the corresponding proof for (pairs of) groups \cite{jolissant}. For the proof, and throughout the rest of the paper, we denote the unitary group of $C(\GG_\max)$ by $\U$. 

\begin{proof}[Proof of Lemma \ref{unitary-lem}]
Denote by $C$ the closed, convex hull of the set 
\[
\Omega=\{\pi(v)\varps(v^*)\xi\mid v\in \U\}.
\]
For any element  $\eta=\sum_{k=1}^n t_k\pi(v_k)\varps(v_k^*)\xi$ in the convex hull of $\Omega$ we have
\begin{align*}
\|\xi-\eta\|&=\|\sum_{k=1}^n t_k(\pi(v_k)\varps(v_k^*)\xi-\xi)\|\\
&\leq \sum_{k=1}^n t_k\|\pi(v_k)\varps(v_k^*)\xi-\xi\|\\
&=\sum_{k=1}^n t_k \|\pi(v_k)\xi-\varps(v_k)\xi\|\leq\delta,
\end{align*}
and hence $\|\xi-\eta\|\leq \delta$ for any $\eta\in C$. Now let $\eta\in C$ be the unique element of minimal norm \cite[Proposition 2.2.1]{KR1}. For every $v\in \U$ we have 
\begin{align*}
\pi(v)\Omega &= \pi(v)\{\pi(u)\varps(u^*)\xi\mid u\in \U\}\\
&=\pi(v)\{\pi(v^*u)\varps(u^*v)\xi\mid u\in \U  \}\\
&=\varps(v)\Omega,
\end{align*}
and hence $\pi(v)C=\varps(v)C$. Since $\pi(v)\eta$ is the element of minimal norm in $\pi(v)C$ and $\varps(v)\eta$ is the ditto element in $\varps(v)C$ we conclude that $\pi(v)\eta=\varps(v)\eta$ for every $v\in \U$. But since the elements in $\U$ span $C(\GG_\max)$ linearly the vector $\eta$ is invariant. 
\end{proof}

\begin{proof}[Proof of Theorem \ref{convergence-thm}]
Assume first that $\hat{\GG}$ has property $\T$ and consider any net $(\varphi_\lambda)_{\lambda\in \Lambda}$ of states on $C(\GG_\max)$ converging pointwise to $\varps$. Denote by  $(H_\lambda,\pi_\lambda,\xi_\lambda)$ the GNS-triple associated with $\varphi_\lambda$. A straight forward calculation reveals that
\begin{align}
|\varphi_\lambda(a)-\varps(a)|^2=\|\pi_\lambda(a)\xi_\lambda-\varps(a)\xi_\lambda\|^2 -(\varphi_\lambda(a^*a)-\varphi_\lambda(a^*)\varphi_\lambda(a))\label{dagger-eq}
\end{align}
for any $a\in C(\GG_\max)$. Note also that the Cauchy-Schwarz inequality implies that $0\leq \varphi_\lambda(a^*a)-\varphi_\lambda(a^*)\varphi_\lambda(a) $ and that this quantity converges to zero. Hence $\lim_\lambda\|\pi_\lambda(a)\xi_\lambda-\varps(a)\xi_\lambda\|=0$ for every $a\in C(\GG_\max)$. Let $\delta>0$ be given. Since $\hat{\GG}$ has property $\T$, Proposition \ref{continuity-prop} allows us to find a Kazhdan pair $(E_0,\delta_0)$ such that any $*$-representation with an $(E_0,\delta_0)$-invariant unit vector $\xi$ has an invariant vector $\eta$ such that $\|\xi-\eta\|\leq \frac{{\delta}}{2}$.
We now claim that the $*$-representation
\[ 	
\pi:=\mathop{\bigoplus}_{\lambda\in \Lambda}{}\pi_{\lambda}  \colon \Pol(\GG)\to B\Bigl(\mathop{\bigoplus}_{\lambda\in \Lambda}^{} H_\lambda\Bigr)
\]
is of this type. To see this, denote by $\tilde{\xi}_{\lambda}$ the image of $\xi_\lambda$ under the natural embedding of $H_\lambda$ into $H:=\oplus_{\mu\in \Lambda}H_\mu$ and note that
\[
\|\pi(a)\tilde{\xi}_\lambda-\varps(a)\tilde{\xi}_\lambda\|=\|\pi_\lambda(a)\xi_\lambda-\varps(a)\xi_\lambda\|\underset{\lambda}{\To}0
\]
for any $a\in C(\GG_\max)$. In particular we get a $\lambda_0\in \Lambda$ such that
\[
\forall \lambda\geq \lambda_0 \ \forall \alpha\in E_0 \ \forall i,j\in \{1,\dots, n_\alpha\}: \|\pi(u_{ij}^\alpha)\tilde{\xi}_\lambda-\varps(u_{ij}^\alpha)\tilde{\xi}_\lambda\|<\delta_0,
\]
and we may therefore find, for each $\lambda\geq\lambda_0$, an invariant vector ${\eta}_\lambda\in H$ such that $\|\tilde{\xi}_\lambda-{\eta}_\lambda\|\leq \frac{{\delta}}{2}$. The equation (\ref{dagger-eq}) now gives
\begin{align*}
|\varphi_\lambda(a)-\varps(a)| & \leq \|\pi_\lambda(a)\xi_\lambda-\varps(a)\xi_\lambda\|\\
&=\|\pi(a)\tilde{\xi}_\lambda-\varps(a)\tilde{\xi}_\lambda\|\\
&=\|\pi(a)(\tilde{\xi}_\lambda-{\eta}_\lambda)+\varps(a)({\eta}_\lambda-\tilde{\xi}_\lambda)\|\\
&\leq \|\pi(a)\|\|\tilde{\xi}_\lambda-{\eta}_\lambda\| +|\varps(a)| \|\tilde{\xi}_\lambda-{\eta}_\lambda\|^2\\
&\leq \delta \|a\|_{\max}
\end{align*}
for every $\lambda\geq \lambda_0$. Hence $(\varphi_\lambda)_{\lambda\in \Lambda}$ converges uniformly to $\varps$ as desired. \\

Assume, conversely, that $\hat{\GG}$ does not have property $\T$ and choose an increasing sequence of subsets $E_n\subseteq \Irred(\GG)$ with union $\Irred(\GG)$. This is possible since $C(\GG)$ is assumed separable so that $\Irred(\GG)$  is a countable set. By Proposition \ref{continuity-prop} we can find $\delta_0>0$ such that for any $n\in \NN$ there exists a Hilbert space $H_n$ and a $*$-representation $\pi_n\colon \Pol(\GG)\to B(H_n)$ which has an $(E_n,\tiny{\frac{1}{n}})$-invariant unit vector $\xi_n$, but such that any invariant vector is a least $\delta_0$ away from $\xi_n$. Define $\varphi_n\colon C(\GG_\max)\to \CC$ by $\varphi_n(a)=\ip{\pi_n(a)\xi_n}{\xi_n}$. Just as above, we get that each $\varphi_n$ satisfies the equation (\ref{dagger-eq}) and by construction of the $E_n$'s it follows that $\lim_n\|\pi_n(a)\xi_n-\varps(a)\xi_n\|=0$ for any $a\in \Pol(\GG)$. Hence $(\varphi_n)_{n\in \NN}$ converges pointwise to $\varps$ on $\Pol(\GG)$ and a standard approximation argument shows that the pointwise convergence then holds on all of $C(\GG_\max)$. Since there are no non-zero invariant vectors within distance $\frac{\delta_0}{2}$ from $\xi_n$, Lemma \ref{unitary-lem} provides us with a  $v_n\in \U$ such that
\[
\|\pi(v_n)\xi_n-\varps(v_n)\xi_n\|> \frac{\delta_0}{2}
\]
Using again the equation (\ref{dagger-eq}) we see that
\[
|\varphi_n(v_n)-\varps(v_n)|^2+(1-|\varphi_n(v_n)|^2)\geq \frac{\delta_0^2}{4},
\]
proving that the convergence can not be uniform.
\end{proof}
The proof of Theorem \ref{convergence-thm} shows that we can get a bit closer to the classical formulation in that we can replace nets with sequences. 

\begin{cor}
The discrete quantum group $\hat{\GG}$ has property $\T$ iff any sequence of states on $C(\GG_\max)$ converging pointwise to the counit converges in the uniform norm.
\end{cor}
\begin{proof}
If $\hat{\GG}$ has property $\T$ the desired conclusion follows from Theorem \ref{convergence-thm}. If, on the other hand, $\hat{\GG}$ does not have property $\T$ the proof of Theorem \ref{convergence-thm} shows how to construct a sequence of states converging pointwise, but not uniformly, to the counit.
\end{proof}

\section{Cocycles and conditionally negative functions}\label{cocycle-section}

The Delorme-Guichardet theorem for groups, stated in the introduction, expresses property $\T$ in terms of vanishing of the first cohomology of the group in question. In order to prove a quantum group version of this result we first introduce the relevant notion of cohomology.

\begin{defi}
Let $\hat{\GG}$ be a discrete quantum group and let $\pi\colon \Pol(\GG)\to B(H)$ be a $*$-representation. A $1$-cocycle for the $*$-representation $\pi$ is a linear map $c\colon \Pol(\GG)\to H$ satisfying
\[
c(ab)=\pi(a)c(b)+c(a)\varps(b),
\]
for all $a,b\in \Pol(\GG)$. A  $1$-cocycle $c$ is called inner if there exists $\xi\in H$ such that $c(a)=\pi(a)\xi-\xi\varps(a)$ for all $a\in \Pol(\GG)$. The set of cocycles $Z^1(\Pol(\GG),H)$ is naturally a complex vector space in which the set of inner cocycles $B^1(\Pol(\GG),H)$ constitutes a subspace, and the first cohomology $H^1(\Pol{\GG},H)$  with coefficients in $H$  is then defined as the space of cocycles modulo the space of inner ones. Finally, a cocycle $c$ is called real if
\[
\ip{c(S(y^*))}{c((Sx)^*)}=\ip{c(x)}{c(y)} \ \text{ for all } x,y\in\Pol(\GG).
\]
\end{defi} 

\begin{rem}
Note that a cocycle $c\colon \Pol(\GG)\to H$ is nothing but a derivation into $H$ where $H$ is considered a $\Pol(\GG)$-bimodule with left action given by $\pi$ and right action given by the counit $\varps$. This is the reason why we from time to time, a bit unconventionally, write the scalar action via $\varps$ on the right. Using the standard description of the first Hochschild cohomology in terms of derivations \cite{loday}, we see that $H^1(\Pol{\GG}, H)$ is exactly the first Hochschild cohomology of $\Pol(\GG)$ with coefficients in the bimodule $_\pi H_\varps$. Throughout the paper, we shall only make use of the \emph{first} Hochschild cohomology group and in the sequel the term cocycle will therefore be used to mean 1-cocycle.
\end{rem}
The following lemma gives an alternative description of the space of inner cocycles and  is a modified version of a result in \cite{peterson-T}.

\begin{lem}\label{bd-is-inner}
If $\pi\colon \Pol(\GG)\to B(H)$ is a $*$-representation and $c\colon \Pol(\GG)\to H$ is a cocycle then $c$ is inner if and only if it is bounded with respect to the norm $\|\cdot\|_\max$ on $C(\GG_\max)$.
\end{lem}
\begin{proof}
First note that both $\pi$ and $\varps$ extend to $C(\GG_\max)$ by definition of the universal norm. It is clear that an inner cocycle is bounded so assume, conversely,  that $c$ extends to $C(\GG_\max)$. We denote the extensions of $\pi,\varps$ and $c$ by the same symbols and define
\[
X=\{c(u)\varps(u^*)\mid u\in \U \},
\]
where $\U$ as before denotes the unitary group of $C(\GG_\max)$. Since $X$ is a bounded set in the Hilbert space $H$ there is a unique Chebyshev center \cite[Lemma 2.2.7]{BHV}; i.e.~there exists a unique $\xi_0\in H$ minimizing the function
\[
H\ni \xi\mapsto \sup\{\|x-\xi\|\mid x\in X\}\in \RR.
\]
Consider now the affine isometric action of  $\mathcal{U}$ on $H$ given by $\alpha(v)(\xi)=\pi(v)\xi +c(v)$. Then for any $v\in \mathcal{U}$ we have that $\alpha(v)\xi_0$ is the Chebyshev center for $\alpha(v)X$ and that $\xi_0\varps(v)$ is the Chebyshev center for $X\varps(v)$. On the other hand
\begin{align*}
\alpha(v)X &= \alpha(v)\{c(u)\varps(u^*)\mid u \in  \mathcal{U}\}\\
&=\alpha(v)\{c(v^*u)\varps(u^*v)\mid u\in \mathcal{U}\}\\
&=\{\pi(v)c(v^*u)\varps(u^*v)+c(v)\mid u\in \mathcal{U}\}\\
&=\{\pi(v)(\pi(v^*)c(u)+c(v^*)\varps(u))\varps(u^*v)+c(v)\mid u\in \mathcal{U}\}\\
&=\{c(u)\varps(u^*)\varps(v)+ \pi(v)c(v^*)\varps(v) + c(v)\mid u\in \mathcal{U}\}\\
&=\{c(u)\varps(u^*)\varps(v)- c(v)\varps(v^*)\varps(v) + c(v)\mid u\in \mathcal{U}\}\\
&=\{c(u)\varps(u^*)\varps(v)\mid u\in \mathcal{U}\}\\
&=X\varps(v),
\end{align*}
and hence $\alpha(v)\xi_0=\xi_0\varps(v)$ for any $v\in \U$. Thus $c(v)=\pi(v)(-\xi_0)-(-\xi_0)\varps(v)$ for any $v\in \mathcal{U}$ and since the elements in $\mathcal{U}$ span $C(\GG_\max)$ linearly we conclude that  $c$ is inner.

\end{proof}
The notion of 1-cocycles on a discrete group $\Gamma$ is intimately linked (see e.g.~section 2.10 in \cite{BHV}) to the notion of conditionally negative definite functions. Recall, that a function $\psi\colon \Gamma\to \RR$ is called conditionally negative definite if $\psi(\gamma)=\psi(\gamma^{-1})$ for every $\gamma\in \Gamma$  and if $\psi$, for any finite subset $\{\gamma_1,\dots,\gamma_n\}\subseteq \Gamma$, furthermore satisfies 
\[
\sum_{i,j=1}^n {\alpha}_i\alpha_j\psi(\gamma_i^{-1}\gamma_j)\leq 0 \text{ for all }  \alpha_1,\dots, \alpha_n\in \RR \text{ with } \sum_{i=1}^n\alpha_i=0. 
\]
The function $\psi$ is said to be normalized if $\psi(e)=0$. Generalizing this to quantum groups we arrive at the following definition.
\begin{defi}
A functional $\psi\colon \Pol(\GG)\to \CC$ is said to be conditionally negative if $\psi(x^*x)\leq 0$ for all $x\in \ker(\varps)$. Moreover, $\psi$ is called normalized if $\psi(1)=0$ and hermitian if $\psi(x^*)=\overline{\psi(x)}$ for all $x\in \Pol(\GG)$.
\end{defi}
The conditionally negative, normalized and hermitian functionals are also cal\-led \emph{infinitesimal generators} because of the following version of Sch{\"o}nberg's theorem.

\begin{thm}[Sch{\"u}rmann, \cite{schurmann-gaussian-states}]
A functional $\psi\colon \Pol(\GG)\to \CC$ is conditionally negative, normalized and hermitian if and only if $\varphi_t=\exp(-t\psi)\colon \Pol(\GG)\to \CC$ is a positive and unital functional for every $t\geq 0$.
\end{thm}
Here positivity of the map $\varphi_t$ simply means that $\varphi_t(x^*x)\geq 0$ for every $x\in \Pol(\GG)$. Note that \cite[Theorem 3.3]{murphy-tuset}  states that such functionals automatically extend to states on $C(\GG_\max)$.  Perhaps the definition of the $\varphi_t$'s require a bit of explanation. For two functionals $\mu,\omega\colon \Pol(\GG)\to\CC$ their convolution product $\omega\star \mu$ is defined as $(\omega\tens\mu)\Delta$. For a single functional $\psi$, the co-semisimplicity of $\Pol(\GG)$ makes the series 
\[
\sum_{k=0}^\infty \frac{(-t)^{k}}{k!}\psi^{\star k}(x)
\]
convergent for each $x\in \Pol(\GG)$ and its sum is denoted $\exp(-t\psi)(x)$. For an infinitesimal generator $\psi$, the family $\varphi_t$ defined above is actually a 1-parameter convolution semigroup of states on $C(\GG_\max)$ converging pointwise to the counit; i.e.~for all $s,t\geq 0$ we have $\varphi_t\star \varphi_s=\varphi_{t+s}$,  $\varphi_0=\varps$  and for every  $x\in \Pol(\GG)$ we have $\varps(x)=\lim_{t\to 0}\varphi_t(x)$.  Such 1-parameter semigroups of states on $C^*$-bialgebras have been studied by Lindsay and Skalski in \cite{skalski-lindsay-convolution} where it is also proved that if 
\[
\lim_{t\to 0}\|\varphi_t-\varps\|=0,
\]
i.e.~if the convergence is uniform,  then the infinitesimal generator $\psi$ is bounded with respect to the norm $\|\cdot\|_\max$. \\

	In the classical case, a group cocycle $c\colon \Gamma \to {_\pi H}$ gives rise to an infinitesimal generator  $\psi\colon \Gamma\to \CC$ by setting $\psi(\gamma)=\|c(\gamma)\|^2$; for quantum groups we have the following analogous result.

\begin{thm}[Vergnioux]\label{vergnioux-thm}
Let $\pi\colon \Pol(\GG)\to B(H)$ be a $*$-representation and $c\colon \Pol(\GG)\to H $ a cocycle. Then  $\psi\colon\Pol(\GG)\to\CC$ defined by 
\[
\psi(x)=\ip{c(x_{(1)})}{c(S(x_{(2)}^*))}
\]
is linear and satisfies 
\begin{align}
\psi(x^*y)=-\ip{c((Sx)^*)}{c(S(y^*))}-\ip{c(y)}{c(x)} \ \text{ for all $x,y\in \ker(\varps)$.}  \label{psi-cond-neg}
\end{align}
If furthermore $c$ is a real cocycle then $\psi$ is an infinitesimal generator; i.e.~$\psi$ is conditionally negative, normalized and hermitian.
\end{thm}	
In the definition of the functional $\psi$ we made use of the so-called Sweedler notation, writing $x_\one\tens x_\two$ for $\Delta x$. We shall use this notation without further elaboration in the following and refer the reader to \cite{klimyk} for a detailed treatment. Theorem \ref{vergnioux-thm} is due to R.~Vergnioux and the author would like to express his gratitude to Vergnioux for communicating it and for allowing its appearance in the present paper. Since the result is not published elsewhere we include Vergnioux's proof, but before doing so a bit of notation is needed. 

\begin{notation}
The dual of the Hilbert space $H$ is denoted $H^\op$ and the inner product in $H$ will be considered both as a sesquilinear form $\ip{\cdot}{\cdot}\colon H\times H\to \CC$ and as a linear map $\ip{\cdot}{\cdot}\colon H\odot H^\op\to \CC$. For $\xi\in H$ we denote by $\xi^\op\in H^\op$ the dual element $\ip{\cdot}{\xi}$ and for $T\in B(H)$ we denote by $T^\op\in B(H^\op)$ the operator $T^\op\xi^\op=(T\xi)^\op$. The symbol $m$ will denote both the multiplication map $\Pol(\GG)\odot \Pol(\GG)\to \Pol(\GG)$ as well as the action $\pi(\Pol(\GG))\odot H\to H$. The antipode in $\Pol(\GG)$ is denoted by $S$ and as usual we denote the counit by $\varps$. Recall that in the general (i.e.~non-Kac) case $S^2\neq \id$, but the relation $S(S(x^*)^*)=x$ always holds. We will often consider the $*$-operation as a self-map of $\Pol(\GG)$ and may therefore write $*(a)$ instead of $a^*$; the above relation involving the antipode may then be written as $S*S\hspace{0.1cm}*=\id$.  Likewise, we will consider $\xi\mapsto \xi^\op$ as a map $\op\colon H\to H^\op$ an write $\op(\xi)$ in stead of $\xi^\op$ whenever convenient. By $\sigma$ we will denote the flip-map on $\Pol(\GG)\odot\Pol(\GG)$ as well as the flip-map on $H\odot H$. Similarly, $\sigma_{(13)}$ will denote the map on a three-fold tensor product which flips the first and the third leg and leaves the middle leg untouched. Throughout this section, we will furthermore make use of the abundance of relations valid in a Hopf $*$-algebra without further reference. These may be found in any standard book on Hopf $*$-algebras; for instance \cite{klimyk}. 
\end{notation}
For the proof of Theorem \ref{vergnioux-thm} we will need a small lemma concerning the interplay between the cocycle and the antipode.

\begin{lem}[Vergnioux]\label{vergnioux-lem}
For $x\in \Pol(\GG)$ with $\varps(x)=0$ we have
\begin{align}
c(Sx) &= -\pi(Sx_{(1)})c(x_\two);\label{c-eq1} \\
c(x) &= -\pi(x_\one)c(Sx_\two);\label{c-eq2}\\
c((Sx)^*)&= - \pi((Sx_\two)^*)c(x_\one^*).\label{c-eq3}
\end{align}
\end{lem}

\begin{proof}
The relation $m(S\tens\id)\Delta x=\varps(x)1$ gives
\begin{align*}
0 &=  c(m(S\tens\id)\Delta x)\\
&=c( (Sx_\one)x_\two)\\
&= \pi(Sx_\one)c(x_\two) +c(Sx_\one)\varps(x_\two)\\
&=\pi(Sx_\one)c(x_\two) + c(S((\id\tens\varps)\Delta x))\\
&=\pi(Sx_\one)c(x_\two) + c(Sx),
\end{align*}
proving equation (\ref{c-eq1}). In the same manner, the equation (\ref{c-eq2}) follows from the formula $m(\id\tens S)\Delta x=\varps(x)1$. The equation (\ref{c-eq3}) follows from (\ref{c-eq2}) and the formula $\Delta S=(S\tens S)\sigma \Delta$:
\begin{align*}
c((Sx)^*) &=-\pi\left(((Sx)^*)_\one\right)c\left(S( ((Sx)^*)_\two )\right)\\
&=\pi\left((Sx_\two)^*\right)c\left(S((Sx_\one)^*)\right)\\
&=\pi((Sx_\two)^*)c(x_\one^*).
\end{align*}
This concludes the proof of the lemma.
\end{proof}

\begin{proof}[Proof of Theorem \ref{vergnioux-thm}]
The functional may be written as $\psi=\ipp{}{} c\tens (\op \ c\ S\ *)\Delta$ which shows that $\psi$ is well defined and linear.
We first prove that $\psi$ satisfies equation (\ref{psi-cond-neg}). Using the cocycle condition we get
\begin{align*}
\psi(x^*y)&=\ipp c\tens (\op\ c\ S\ *)\Delta(x^*y)\\
&= \ip{c(x_\one^*y_\one)}{c((Sx_\two)S(y_\two^*))}\\
&=\ip{\pi(x_\one^*)c(y_\one)}{\pi(Sx_\two)c(S(y_\two^*))}+\ip{\pi(x_\one^*)c(y_\one)}{c(Sx_\two)\varps(S(y_\two^*))}\\
&+ \ip{c(x_\one^*)\varps(y_\one)}{\pi(Sx_\two)c(S(y_\two^*))}+\ip{c(x_\one^*)\varps(y_\one)}{c(Sx_\two)\varps(S(y_\two^*))}.
\end{align*}
We now treat the four terms one by one.  Since $\varps(x)=0$, the first term vanishes:
\begin{align*}
\ip{\pi(x_\one^*)c(y_\one)}{\pi(Sx_\two)c(S(y_\two^*))} &= \ip{c(y_\one)}{\pi(m(\id\tens S)\Delta x)c(S(y_\two^*))}\\
&=\ip{c(y_\one)}{\pi(\varps(x)1)c(S(y_\two^*))}\\
&=0.
\end{align*}
Using the formula (\ref{c-eq1}) and the fact that $\varps S=\varps $, the second term becomes
\begin{align*}
\ip{\pi(x_\one^*)c(y_\one)}{c(Sx_\two)\varps(S(y_\two^*))} &= \ip{\varps(y_\two)c(y_\one)}{\pi(x_\one)c(Sx_\two)}\\
&=-\ip{c((\id\tens\varps)\Delta y) }{c(x)}\\
&=-\ip{c(y)}{c(x)}.
\end{align*}
Similarly, using the formula (\ref{c-eq3}) the third term becomes
\begin{align*}
\ip{c(x_\one^*)\varps(y_\one)}{\pi(Sx_\two)c(S(y_\two^*))} &= \ip{\pi((Sx_\two)^*)c(x_\one^*)}{\varps(y_\one^*)c(S(y_\two^*))}\\
&=-\ip{c((Sx)^*)}{c(S((\varps\tens\id)\Delta (y^*)))}\\
&=-\ip{c((Sx)^*)}{c(S(y^*))}.
\end{align*}
We are therefore done if we can show that the fourth term vanishes. This follows from the assumption $\varps(y)=0$ and the following calculation:
\begin{align*}
\ip{c(x_\one^*)\varps(y_\one)}{c(Sx_\two)\varps(S(y_\two^*)))} &= \varps(y_\one)\varps(y_\two)\ip{c(x_\one^*)}{c(Sx_\two)}\\
&=\varps((\id\tens\varps)\Delta y)\ip{c(x_\one^*)}{c(Sx_\two)}\\
&=0.
\end{align*}
Hence $\psi$ satisfies  (\ref{psi-cond-neg}). Assume now that $c$ is a real cocycle. The equation (\ref{psi-cond-neg}) then gives that $\psi(x^*x)=-2\|c(x)\|^2$ whenever $x\in \ker(\varps)$ which shows that $\psi$ is conditionally negative. Since $c(1)=0$ it is clear that $\psi$ is normalized. That $\psi$ is hermitian is seen by the following calculation.
\begin{align*}
\overline{\psi(x^*)} &= \overline{\ipp (c\tens (\op \ c \ S \ *)) x_\one^*\tens x_\two^*}\\
&=\overline{\ip{c(x_\one^*)}{c(Sx_\two)}}\\
&=\ip{c(Sx_\two)}{c(x_\one^*)}\\
&=-\ip{\pi(S({x_\two}_\one))c({x_\two}_\two)}{c(x_\one^*)}\tag{by (\ref{c-eq1})}\\
&=-\ip{c({x_\two}_\two)}{\pi(S({x_\two}_\one)^*)  c(x_\one^*)}\\
&=-\ipp c({x_\two}_\two)\tens (\pi((S{x_\two}_\one)^*)^\op c(x_\one^*)^\op\\
&=-\ipp c\tens [m((\op \ \pi \hspace{0.05cm} * \hspace{0.05cm} S)\tens (\op \ c \ *))]  ({x_\two}_\two\tens {x_\two}_\one\tens x_\one)\\
&=-\ipp c\tens [m((\op \ \pi \hspace{0.05cm} *\hspace{0.05cm} S)\tens (\op \ c \ *))]\sigma_{(13)}  (x_\one\tens {x_\two}_\one\tens {x_\two}_\two)\\
&=-\ipp c\tens [m((\op \ \pi \hspace{0.05cm}*\hspace{0.05cm} S)\tens (\op \ c \ *))]\sigma_{(13)}  ({x_\one}_\one \tens {x_\one}_\two\tens {x_\two})\\
&=-\ip{c(x_\two)}{\pi((S({x_\one}_\two))^*)c({x_\one}_\one^*)}\\
&=\ip{c(x_\two)}{c((Sx_\one)^*)}\tag{by (\ref{c-eq3})}\\
&=\ip{c(S(S(x_\two^*)^*))}{c((Sx_\one)^*)}\\
&=\ip{c(x_\one)}{c(S(x_\two^*))}\tag{$c$ real}\\
&=\psi(x).
\end{align*}
This concludes the proof of Theorem \ref{vergnioux-thm}.
\end{proof}

\section{The Delorme-Guichardet Theorem}
In this section we prove our main result which characterizes property $\T$ of a discrete quantum group in terms of its first cohomology groups.
\begin{thm}\label{delorme-guichardet}
For a discrete quantum group $\hat{\GG}$ the following are equivalent.
\begin{itemize}
\item[(i)] $\hat{\GG}$ has property $\T$.
\item[(ii)] $\hat{\GG}$ is Kac and every normalized, hermitian, conditionally negative functional $\psi\colon \Pol(\GG)\to\CC$ is bounded with respect to $\|\cdot\|_\max$.
\item[(iii)] For every $*$-representation $\pi\colon \Pol(\GG)\to B(H)$ the first cohomology group $H^1(\Pol(\GG),H)$ vanishes.
\end{itemize}
\end{thm}
Before giving the proof we introduce a topology on the space of 1-cocycles. Let $\pi\colon\Pol(\GG)\to B(H)$ be a $*$-representation and
define, for each finite subset $E\subseteq \Irred(\GG)$,  a seminorm on $Z^1(\Pol(\GG),H)$ by
\[
\|c\|_E=\sup\{\|c(u_{ij}^\alpha)\|\mid \alpha\in E, 1\leq i,j\leq n_\alpha\}.
\]
Since the matrix coefficients span $\Pol(\GG)$ linearly, it is a routine to check  that $Z^1(\Pol(\GG),H)$ becomes a Frechet space when endowed with the topology arising from this family of seminorms.  This topology captures the existence of almost invariant vectors for $\pi$ by means of the following lemma, which  generalizes a result of Guichardet \cite[Th{\'e}or{\`{e}}me 1]{guichardet}.
\begin{lem}\label{frechet-lem}
Assume that $\pi$ does not have non-zero invariant vectors. Then $B^1(\Pol(\GG),H)$ is closed in $Z^1(\Pol(\GG),H)$ if and only if $\pi$ does not have almost invariant vectors. 
\end{lem}
\begin{proof}
Consider the map $\Phi\colon H\to B^1(\Pol(\GG),H)$ mapping a vector $\xi$ to the corresponding inner cocycle. Then $\Phi$ is linear, continuous and surjective and since $\pi$ is assumed to have no non-zero fixed vectors it follows that $\Phi$ is also injective. Assume first that $\pi$ does not have almost invariant vectors either. Then there exists $E_0\subseteq \Irred(\GG)$ and $\delta_0>0$ such that 
\begin{align}
\|\Phi(\xi)\|_{E_0}=\sup\{\|\pi(u_{ij}^\alpha)\xi-\varps(u_{ij}^\alpha)\xi\| \mid \alpha\in E_0, 1\leq i,j\leq n_\alpha \}      \geq \delta_0 \|\xi\| \label{E0-ineq} 
\end{align}
for all $\xi\in H$. Let $c$ be in the closure of $B^1(\Pol(\GG),H)$ and choose a sequence $(\xi_n)_{n\in \NN}$ such that $\Phi(\xi_n)$ converges to $c$ in the Frechet topology. Then, in particular, $(\Phi(\xi_n))_{n\in \NN}$ is a Cauchy sequence with respect to the seminorm $\|\cdot\|_{E_0}$ and  by (\ref{E0-ineq}) this implies that $(\xi_n)_{n\in \NN}$ is a Cauchy sequence in $H$;  hence it has a limit point $\xi\in H$. By continuity of $\Phi$ we conclude that $c=\Phi(\xi)$ and therefore $c$ is inner and $B^1(\Pol(\GG),H)$ closed. Conversely, assume that $B^1(\Pol(\GG),H)$ is closed and therefore a sub-Frechet space in $Z^1(\Pol(\GG),H)$. Then the open mapping theorem \cite[Corollaries 2.12]{rudin-functional-analysis} implies that $\Phi\colon H\to B^1(\Pol(\GG),H)$ is bi-continuous and thus bounded away from zero in at least  one of the seminorms; say  $\|\cdot\|_{E_0}$.  
Hence there exists a $\delta_0>0$ such that
\[
\sup\{\|\pi(u_{ij}^\alpha)\xi-\varps(u_{ij}^\alpha)\xi\|\mid \alpha\in E_0,1\leq i,j\leq n_\alpha\}=\|\Phi(\xi)\|_{E_0}\geq \delta_0 \|\xi\| 
\]
for all $\xi\in H$, and therefore  $\pi$ can not have almost invariant vectors.
\end{proof}
We are now ready to give the proof of the Delorme-Guichardet theorem.
\begin{proof}[Proof of Theorem \ref{delorme-guichardet}]
We first prove (i)$\Rightarrow$(ii). Assume therefore that $\hat{\GG}$ has property $\T$ and let $\psi\colon \Pol(\GG)\to \CC$ be normalized, conditionally negative and hermitian. By exponentiation (see Section \ref{cocycle-section}) we get a 1-parameter family of states on $C(\GG_\max)$ converging pointwise to the counit, and by Theorem \ref{convergence-thm} the convergence has to be uniform. Applying \cite[Proposition 2.3]{skalski-lindsay-convolution}, this implies that the infinitesimal generator $\psi$ is bounded. That $\hat{\GG}$ is Kac follows from Theorem \ref{fima-thm}. \\

Next we prove (ii)$\Rightarrow$(iii). Assume therefore that $\hat{\GG}$ is Kac and that every infinitesimal generator is bounded, and let furthermore a $*$-representation $\pi\colon\Pol(\GG)\to B(H)$ as well as a cocycle $c\colon \Pol(\GG)\to H$ be given. By Lemma \ref{bd-is-inner} it suffices to show that $c$ is bounded with respect to the norm $\|\cdot\|_\max$. Since $\hat{\GG}$  is assumed Kac the antipode $S\colon \Pol(\GG)\to \Pol(\GG)$ is $*$-preserving. The $*$-representation $\pi$ therefore gives rise to a dual $*$-representation 
$\pi^\op\colon \Pol(\GG)\to B(H^\op)$ on the dual Hilbert space $H^\op$ given by $\pi^\op(a)\xi^\op=(\pi(Sa^*)\xi)^\op$ and the map $c^\op\colon\Pol(\GG)\to H^\op$ given by $c^\op(a)=(c(Sa^*))^\op$ is a cocycle for this representation. Moreover, it is easy to check that 
\[
\Pol(\GG)\ni a \longmapsto (c(a),c^\op(a))\in H\oplus H^\op
\]
is a real $(\pi\oplus \pi^\op)$-cocycle which is bounded if and only if $c$ is bounded. It therefore suffices to treat the case where the cocycle $c$ is real. In this case, Theorem \ref{vergnioux-thm} provides us with a conditionally negative functional $\psi\colon \Pol(\GG)\to \CC$ such that
$
\psi(x^*x)= -2\|c(x)\|^2
$
for all $x\in \ker(\varps)$. By assumption the functional $\psi$ is bounded with respect to $\|\cdot\|_\max$ so for $x\in \ker(\varps)$ we get
\[
\|x\|_\max^2\|\psi\|\geq |\psi(x^*x)|=2\|c(x)\|^2,
\]
and hence the restriction $c_0$ of $c$ to $\ker(\varps)$ extends boundedly to a map $\tilde{c}_0$ on the closure $J$ of $\ker(\varps)$ inside $C(\GG_\max)$. The ideal $J$ is exactly the kernel of the extension of the counit $\varps\colon C(\GG_\max)\to \CC$ and therefore the map
\[
C(\GG_\max)\ni a \mapsto \tilde{c}_0(a-\varps(a)1)\in H
\]
is bounded and extends $c$.\\

To prove (iii)$\Rightarrow$(i), assume that all the first cohomology groups vanish and let $\pi\colon \Pol(\GG)\to B(H)$ be a representation without non-zero invariant vectors. By assumption $H^1(\Pol(\GG),H)$ vanishes so in particular $B^1(\Pol(\GG),H)$ is closed in the Frechet topology on $Z^1(\Pol(\GG),H)$, and by Lemma \ref{frechet-lem} this implies that $\pi$ does not have almost invariant vectors; thus $\hat{\GG}$ has property $\T$.

\end{proof}

\section{\texorpdfstring{An application to $L^2$-invariants}{An application to L2-invariants}}\label{L2-section}
In this section we prove that the first $L^2$-Betti number of a discrete quantum group  with 
property $\T$ vanishes. The corresponding statement for groups  was known to Gromov \cite{gromov-asymptotic-invariants}, but the first detailed proof was given by Bekka and Valette in  \cite{bekka-valette-group-cohomology}. The modern homological algebraic approach to $L^2$-invariants developed by L{\"u}ck \cite{Luck02} (see also \cite{thom}, \cite{peterson-thom}, \cite{thom-rank-metric}, \cite{sauer-betti-of-groupoids},\cite{farber}) provides a different proof of this result; it can, for instance, easily be deduced from \cite[Theorem 2.2]{peterson-thom}. In this section we show how this argument can be adapted to the quantum group context. Before doing so, we briefly remind the reader of the necessary definitions concerning $L^2$-Betti numbers for groups and quantum groups. \\

	For a discrete group $\Gamma$, its $L^2$-Betti numbers can be described/defined in purely algebraic terms \cite{Luck02} as $\bet_p(\Gamma)=\dim_{\L(\Gamma)}\Tor_p^{\CC\Gamma}(\L(\Gamma),\CC)$ where $\dim_{\L(\Gamma)}(-)$ is L{\"u}ck's extended Murray-von Neumann dimension. For a discrete quantum group $\hat{\GG}$ of Kac type it is therefore natural to define its $L^2$-Betti numbers as
\[
\bet_p(\hat{\GG})=\dim_{L^\infty(\GG)}\Tor_p^{\Pol(\GG)}(L^\infty(\GG),\CC),
\] 
where $\dim_{L^\infty(\GG)}(-)$ is the Murray-von Neumann dimension arising from the tracial Haar state. These $L^2$-Betti numbers have been studied in \cite{thom-collins}, \cite{vergnioux-paths-in-cayley}, \cite{quantum-betti} and \cite{coamenable-betti}\footnote{Notice that $\bet_p(\hat{\GG})$ is denoted $\bet_p(\GG)$ in \cite{quantum-betti} and \cite{coamenable-betti}.},  and the aim of the present section is to prove the following.

\begin{cor}\label{beta-et=0}
If $\hat{\GG}$ has property $\T$ then $\bet_1(\hat{\GG})=0$.
\end{cor}
Note that if $\hat{\GG}$ has property $\T$ then $\GG$ is automatically of Kac type, and its Haar state $h\colon L^\infty(\GG)\to\CC$ is therefore a trace, so that the first $L^2$-Betti number is defined. As mentioned above, the proof of Corollary \ref{beta-et=0} follows the lines of the corresponding proof in \cite{peterson-thom}. During the proof we will have to consider dimensions of both right and left modules for $L^\infty(\GG)$, and to avoid confusion we will let $\dim_{L^\infty(\GG)}(X)$ denote the dimension of a left module $X$ whereas $\dim_{L^\infty(\GG)^\op}(Y)$ will denote the dimension of a right module $Y$.

\begin{proof}[Proof of Corollary \ref{beta-et=0}]
Denote by $M(\GG)$ the $*$-algebra of closed, densely defined (potentially unbounded) operators affiliated with 
$L^\infty(\GG)$. This is a self-injective and von Neumann regular ring and tensoring $L^\infty(\GG)$-modules with $M(\GG)$ is a flat and dimension preserving functor \cite{reich01}.  Therefore
\begin{align*}
\bet_1(\hat{\GG})&=\dim_{L^\infty(\GG)} M(\GG)\odot_{L^\infty(\GG)}\Tor_1^{\Pol(\GG)}(L^\infty(\GG),\CC)\\
&=\dim_{L^\infty(\GG)}\Tor_1^{\Pol(\GG)}(M(\GG),\CC).
\end{align*}
By \cite[Corollary 3.4]{thom}  we have
\[
\dim_{L^\infty(\GG)}\Tor_1^{\Pol(\GG)}(M(\GG),\CC)= \dim_{L^\infty(\GG)^\op}\Hom_{M(\GG)}(\Tor_1^{\Pol(\GG)}(M(\GG),\CC),M(\GG)),
\]
and using the self-injectiveness of $M(\GG)$ (see e.g.~\cite[Theorem 3.5]{thom}  and its proof) we get an isomorphism of right $M(\GG)$-modules
\[
\Hom_{M(\GG)}(\Tor_1^{\Pol(\GG)}(M(\GG),\CC),M(\GG))\simeq \Ext^1_{\Pol(\GG)}(\CC,M(\GG)).
\]
By considering the bar-resolution of the trivial $\Pol(\GG)$-module $\CC$, one sees that $\Ext^1_{\Pol(\GG)}(\CC,M(\GG))$ may also be computed as the first Hochschild cohomology  $H^1(\Pol(\GG),M(\GG))$ where $M(\GG)$ carries the natural left action of $\Pol(\GG)$ and right action given by the counit. We therefore arrive at the formula
\[
\bet_1(\hat{\GG})=\dim_{L^\infty(\GG)^\op}H^1(\Pol(\GG), M(\GG)).
\]
Since $\hat{\GG}$ has property $\T$, Theorem \ref{delorme-guichardet} implies that {$H^1(\Pol(\GG),L^2(\GG))$} vanishes and we are therefore done if we can prove that
\begin{align*}
\dim_{L^\infty(\GG)^\op}H^1(\Pol(\GG),L^2(\GG)) &=\dim_{L^\infty(\GG)^\op}H^1(\Pol(\GG),M(\GG)).
\end{align*}
To see this we consider the following diagram of right $L^\infty(\GG)$-modules:
\[
\xymatrix@=0.55cm{
0 \ar[r] & B^1(\Pol(\GG), L^\infty(\GG)) \ar[r] \ar[d] & Z^1(\Pol(\GG), L^\infty(\GG)) \ar[d] \ar[r] & H^1(\Pol(\GG), L^\infty(\GG))\ar[r] \ar[d] & 0\\
0 \ar[r] & B^1(\Pol(\GG), L^2(\GG)) \ar[r] \ar[d] & Z^1(\Pol(\GG), L^2(\GG)) \ar[d] \ar[r] & H^1(\Pol(\GG), L^2(\GG))\ar[r] \ar[d] & 0\\
0 \ar[r] & B^1(\Pol(\GG), M(\GG)) \ar[r]  & Z^1(\Pol(\GG), M(\GG))  \ar[r] & H^1(\Pol(\GG), M(\GG))\ar[r]  & 0
}
\]
The rows in this diagram are exact by definition and the two first columns clearly consist of inclusions. We now prove that
\begin{align}
\dim_{L^\infty(\GG)^\op}B^1(\Pol(\GG), L^\infty(\GG))&= \dim_{L^\infty(\GG)^\op} B^1(\Pol(\GG), M(\GG));\label{eq1}  \\
\dim_{L^\infty(\GG)^\op}Z^1(\Pol(\GG), L^\infty(\GG))&= \dim_{L^\infty(\GG)^\op} Z^1(\Pol(\GG), M(\GG))\label{eq2},
\end{align}
and the result then follows from additivity (\cite[Theorem 6.7]{Luck02}) of the dimension function $\dim_{L^\infty(\GG)^\op}(-)$. To prove the equality (\ref{eq1}), notice that the first column identifies with the inclusions $L^\infty(\GG)\subseteq L^2(\GG)\subseteq M(\GG)$ so it suffices to see that 
\[
\dim_{L^\infty(\GG)^\op}\bigl(M(\GG)/L^\infty(\GG)\bigr)=0.
\] 
By  \cite[Theorem 2.4]{sauer-betti-of-groupoids}, it is enough to see that	for every $\xi \in M(\GG)$ and every $\delta>0$ there exists a projection $p\in L^\infty(\GG)$ such that $h(p)\geq 1-\delta$ and $\xi p\in L^\infty(\GG)$. But this follows from the fact all the spectral projections of the absolute value of $\xi$ are in $L^\infty(\GG)$.   
To prove the equality (\ref{eq2}), consider a cocycle $c\colon \Pol(\GG)\to M(\GG)$ and a $\delta>0$. Again by \cite[Theorem 2.4]{sauer-betti-of-groupoids}, we have to find a projection $p\in L^\infty(\GG)$ such that $h(p)\geq 1-\delta$ and $c(-)p\in Z^1(\Pol(\GG),L^\infty(\GG))$. For this, consider the set of matrix coefficients 
\[
\{u_{ij}^\alpha \mid \alpha\in I,1\leq i,j\leq n_\alpha\}.
\]
Since $C(\GG)$ is assumed separable, this set is at most countable so we may choose a sequence of numbers $\delta_{ij}^\alpha>0$ such that 
\[
\sum_{\alpha\in I}\sum_{i,j=1}^{n_\alpha}\delta_{ij}^\alpha\leq\delta.
\]
For each $u_{ij}^\alpha$ we have $c(u_{ij}^\alpha)\in M(\GG)$ and hence we can find a projection $p_{ij}^\alpha\in L^\infty(\GG)$ such that $c(u_{ij}^\alpha)p_{ij}^\alpha\in L^\infty(\GG)$ and such that $h(p_{ij}^\alpha)\geq 1-\delta_{ij}^\alpha$. Let $p$ be the infimum of all these projections. We then have
\[
h(1-p)=h\Bigl(1- \bigwedge_{\alpha,i,j}p_{ij}^\alpha\Bigr)=h\Bigl(\bigvee_{\alpha,i,j}1-p_{ij}^\alpha\Bigr)\leq \sum_{\alpha,i,j}h(1-p_{ij}^\alpha)\leq\delta.
\]
Since the set $\{u_{ij}^\alpha \mid \alpha\in I,1\leq i,j\leq n_\alpha\}$ spans $\Pol(\GG)$ linearly we also have that $c(-)p$ is a cocycle with values in $L^\infty(\GG)$ and the proof is complete.
\end{proof}

Turning things around, vanishing of the first $L^2$-Betti number may be turned into an honest vanishing of cohomology result.

\begin{cor}\label{H=0-cor}
If $\hat{\GG}$ is non-amenable and of Kac type with $\bet_1(\hat{\GG})=0$ then $H^1(\Pol(\GG),L^2(\GG))$ vanishes.
\end{cor}
Recall that a discrete quantum group $\hat{\GG}$ is said to be amenable if the counit $\varps\colon \Pol(\GG)\to \CC$ extends to a bounded character on $C(\GG_\red)$. A detailed study of this notion may be found in \cite{murphy-tuset} and \cite{tomatsu-amenable}. Note also that the results in \cite{vergnioux-paths-in-cayley} show that Corollary \ref{H=0-cor} applies to the duals of the free orthogonal quantum groups $O_n^+$ for $n\geq 3$. The proof Corollary \ref{H=0-cor} is again a modification of the corresponding proof in \cite{peterson-thom}.
\begin{proof}
Since $\bet_1(\hat{\GG})=0$ we have\footnote{See the beginning of the proof of Corollary \ref{beta-et=0}.} that
\[
\dim_{L^\infty(\GG)} \Tor_1^{\Pol(\GG)}(M(\GG),\CC)=0,
\]
and by  \cite[Corollary 3.3]{thom} this implies vanishing of the dual $M(\GG)$-module 
\[
\Hom_{M(\GG)}(\Tor_1^{\Pol(\GG)}(M(\GG),\CC),M(\GG)).
\]
As in the proof of Corollary \ref{beta-et=0} we have an isomorphism of right $M(\GG)$-modules
\[
\Hom_{M(\GG)}(\Tor_1^{\Pol(\GG)}(M(\GG),\CC),M(\GG))\simeq H^1(\Pol(\GG),M(\GG)),
\]
and hence every cocycle with values in $M(\GG)$ is inner. Denote by $J\in B(L^2(\GG))$ the modular conjugation arising from the tracial state $h$ and recall that for each $x\in L^2(\GG)$ the operator $L(x)^0\colon \Lambda a\mapsto J\lambda(a)^*Jx$ is pre-closed and its closure $L(x)$ is affiliated with $L^\infty(\GG)$. We therefore have $L^2(\GG)$ embedded into $M(\GG)$ via the map $L\colon x\mapsto L(x)$ which is easily seen to be an embedding of $\Pol(\GG)$-bimodules. For a given cocycle $c\colon \Pol(\GG)\to L^2(\GG)$ we can therefore find an affiliated operator $\xi\in M(\GG)$ such that $L(c(a))= \lambda(a)\xi-\xi\varps(a)$ for every $a\in \Pol(\GG)$.  Choose now an increasing sequence of projections $p_n\in L^\infty(\GG)$ such that $\xi p_n \in L^\infty(\GG)$ for every $n\in \NN$ and such that $(p_n)_{n\in \NN}$ converges in the strong operator topology to 1. For each $n\in\NN$ and $a\in \Pol(\GG)$ we now have
\[
L(c(a))p_n=\lambda(a)(\xi p_n)-(\xi p_n)\varps(a),
\]
and evaluating the operators on $\Lambda(1)$ we get
\[
Jp_n J(c(a))= \lambda(a)\Lambda(\xi p_n)-\Lambda(\xi p_n)\varps(a).
\]
But also $(Jp_n J)_{n\in \NN}$ converges in the strong operator topology to 1 and hence
\[
c(a)=\lim_{n\to\infty} \bigl( \lambda(a)\Lambda(\xi p_n)-\Lambda(\xi p_n)\varps(a)\bigr),
\]
which proves that $c$ is the pointwise limit of a sequence of inner cocycles with values in $L^2(\GG)$. But since $\hat{\GG}$ is non-amenable the left regular representation $\lambda$ can not
have almost invariant vectors and by Lemma \ref{frechet-lem} the space $B^1(\Pol(\GG),L^2(\GG))$ is therefore closed in the Frechet topology on $Z^1(\Pol(\GG),L^2(\GG))$. This topology is exactly the topology of pointwise convergence and we conclude that $c\in Z^1(\Pol(\GG),L^2(\GG))$ is inner.
\end{proof}
\begin{rem}
For further quantum group applications of property $\T$ we refer the reader to \cite{exotic-norms} where a connection with Bekka's notion of property $\T$ for $C^*$-algebras \cite{bekka-T-for-C} is established and used to construct new examples  of completions of $\Pol(\GG)$ that result in $C^*$-algebraic quantum groups.
\end{rem}

\def\cprime{$'$}

\end{document}